\documentclass{article}
\usepackage{amssymb}
\usepackage{amsmath}
\usepackage{amsfonts}
\usepackage{graphicx}
\usepackage{graphics}

\newtheorem{theorem}{Theorem}[section]
\newtheorem{corollary}[theorem]{Corollary}
\newtheorem{definition}[theorem]{Definition}

\newtheorem{lemma}[theorem]{Lemma}
\newtheorem{notation}[theorem]{Notation}
\newtheorem{proposition}[theorem]{Proposition}

\newfont{\zb}{cmssi12 scaled 1000}

\def\Pr{{\bf Proof:\ }}
\def\mod{{\rm mod\ }}
\def\dom{{\rm dom\ }}
\def\im{{\rm im\ }}
\def\rank{{\rm rank\ }}
\def\fbx{\hspace*{4cm}\hfill${}^{\rule{2mm}{2mm}}$}

\title{On the semigroup of all partial fence-preserving injections on a finite set}

\author{Ilinka Dimitrova \\
\textit{Faculty of Mathematics and Natural Science}\\
\textit{South-West University "Neofit Rilski"}\\
\textit{2700 Blagoevgrad, Bulgaria}\\
\textit{email: ilinka\_dimitrova@swu.bg}\\
~~\\
J\"{o}rg Koppitz \footnote{Supported
by the Alexander von Humboldt Foundation}\\
\textit{Institute of Mathematics}\\
\textit{University of Potsdam}\\
\textit{14476 Potsdam, Germany}\\
\textit{email: koppitz@uni-potsdam.de}}

\begin{document}

\maketitle

\begin{abstract}
For $n \in \mathbb{N}$, let $X_n = \{a_1, a_2, \ldots, a_n\}$ be an $n$ - element set
and let $\textbf{F} = (X_n; <_f)$ be a fence, also called a zigzag poset.
As usual, we denote by $I_n$ the symmetric inverse semigroup on $X_n$.
We say that a transformation $\alpha \in I_n$ is
\textit{fence-preserving} if $x <_f y$ implies that $x\alpha <_f y\alpha$, for all $x, y$ in the domain of $\alpha$.
In this paper, we study the semigroup $PFI_n$ of all partial fence-preserving injections of $X_n$ and its subsemigroup
$IF_n = \{\alpha \in PFI_n : \alpha^{-1}\in PFI_n\}$.
Clearly, $IF_n$ is an inverse semigroup and contains all regular elements of $PFI_n.$
We characterize the Green's relations for the semigroup $IF_n$.
Further, we prove that the semigroup $IF_n$ is generated by its elements with \rank$\geq n-2$.
Moreover, for $n \in 2\mathbb{N}$ we find the least generating set and calculate the rank of $IF_n$.\\
\end{abstract}

\textit{Keywords:} finite transformation semigroup, fence-preserving transformations,
inverse semigroup, Green's relations, generators, rank \\

2010 Mathematics Subject Classification: 20M20, 20M18
\\

\section{Introduction and Preliminaries}

For $n \in \mathbb{N}$, let $X_n = \{a_1, a_2, \ldots, a_n\}$ be an $n$ - element set.
As usual, we denote by $I_n$ the symmetric inverse semigroup on $X_n$, i.e. the partial one-to-one transformation
semigroup on $X_n$ under composition of mappings. The importance of $I_{n}$ to inverse semigroup theory may
be likened to that of the symmetric group $S_{n}$ to group theory. Every finite
inverse semigroup $S$ is embeddable in $I_{n}$, the analogue of Cayley's theorem
for finite groups, and to the regular representation of finite semigroups.
Thus, just as the study of symmetric, alternating and dihedral groups has
made a significant contribution to group theory, so has the study of various
subsemigroups of $I_{n}$, see for example \cite{BP,DK,Fern,GM,Umar}.

Let $\textbf{F} = (X_n; <_f)$ be a \textit{fence}, also called a \textit{zigzag poset}, i.e. a partially ordered set
in which the order relation forms a path with alternating orientations:
$$
a_1 <_f a_2 >_f a_3 <_f \cdots a_n
$$
or
$$
a_1 >_f a_2 <_f a_3 >_f \ldots a_n.
$$
Every element of $\textbf{F}$ is either maximal or minimal.
A fence $\textbf{F}$ is called an up-fence (respectively a down-fence) if $a_1 <_f a_2$ (respectively $a_1 >_f a_2$).
In this paper, without loss of generality, we consider an up-fence.
\begin{center}
\includegraphics[scale=.47]{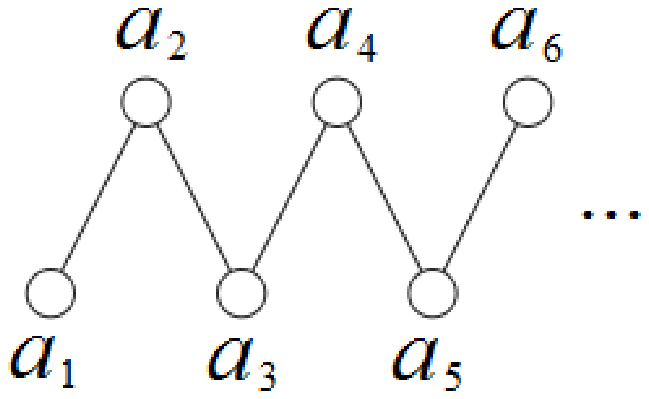}
\end{center}

Several authors have investigated the number of order-preserving maps from fences to themselves, or to fences of other sizes, see for example \cite{CV,Farl,Rutk1,Rutk2}.
Recently, regular semigroups of transformations preserving a fence were characterized in \cite{JS,TSC}.

We begin by recalling some notations and definitions that will be
used in the paper. For standard concepts in semigroup and symmetric inverse semigroup
theory, see for example \cite{How} and \cite{Limp}. We denote by $\dom\alpha $ and $\im\alpha $ the domain
and the image (range) of $\alpha \in I_{n}$, respectively. The natural number
$\rank\alpha :=\left\vert \im\alpha \right\vert$ is called the rank
of $\alpha$. The inverse element of $\alpha$ is denoted by $\alpha^{-1}$.
For a subset $Y \subseteq X_n$, we denote by $id|_Y$ the identity mapping on $Y$.
Clearly, if $Y = X_n$ then $id|_{X_n} =: id$ is the identity mapping on $X_n$.
For a subset $A\subseteq I_{n}$, we denote by $\left\langle A\right\rangle $ the subsemigroup of $I_{n}$
generated by $A$.
We say that a transformation $\alpha \in I_n$ is
\textit{fence-preserving} if $x <_f y$ implies that $x\alpha <_f y\alpha$, for all $x, y \in \dom \alpha$.
We denote by $PFI_n$ the subsemigroup of $I_n$ of all partial fence-preserving injections of $X_n$.
Note that the semigroup $PFI_n$ is not inverse. For example
$$\alpha = \left(
             \begin{array}{ccccc}
               1 & 2 & 4 & 5 & 6 \\
               3 & 2 & 6 & 5 & 4 \\
             \end{array}
           \right) \in PFI_6, ~~\mbox{but}~~
\alpha^{-1} = \left(
             \begin{array}{ccccc}
               2 & 3 & 4 & 5 & 6 \\
               2 & 1 & 6 & 5 & 4 \\
             \end{array}
           \right) \notin PFI_6.
$$
Let $IF_n$ be the set of all $\alpha \in PFI_n$ such that $\alpha^{-1}\in PFI_n$.
Clearly, $IF_n$ is the set of all $\alpha \in PFI_n$ with $x <_f y$ if and only if $x\alpha <_f y\alpha$, for all $x, y \in \dom \alpha$.
Hence, $IF_n$ is an inverse subsemigroup of $PFI_n$ and contains all regular elements of $PFI_n$.
In section 2, we characterize the Green's relations for the inverse semigroup $IF_n$.
Further, we prove that the semigroup $IF_n$ is generated by its elements with \rank$\geq n-2$.
Moreover, for $n \in 2\mathbb{N}$ we find the least generating set and calculate the rank of $IF_n$.

\section{Green's Relations}

In this section, we characterize the Green's relations $\mathcal{R}$, $\mathcal{L}$, $\mathcal{H}$, and $\mathcal{J}$ on $IF_n$.
Since $IF_n$ is an inverse subsemigroup of $I_n$, for $\alpha ,\beta \in IF_n$, it
holds:

\begin{enumerate}
\item $\alpha \mathcal{L}\beta $ if and only if $\im \alpha =\im \beta $.

\item $\alpha \mathcal{R}\beta $ if and only if $\dom \alpha =\dom \beta $.

\item $\alpha\mathcal{H}\beta$ if and only if $\dom \alpha =\dom \beta$ and $\im \alpha =\im \beta$.
\end{enumerate}

It remains to describe the relation $\mathcal{J}$, since this relation is different for the semigroups $I_n$ and $IF_n$.
For example, let
$$\alpha = \left(
             \begin{array}{cccc}
               1 & 4 & 5 & 6 \\
               2 & 6 & 5 & 4 \\
             \end{array}
           \right), ~~
\beta = \left(
             \begin{array}{cccc}
               1 & 2 & 5 & 6 \\
               5 & 6 & 1 & 2 \\
             \end{array}
           \right) \in IF_6.
$$
Then $\rank \alpha = \rank \beta$, but $\alpha$ and $\beta$ are not $\cal J$ related.

\begin{definition}
\rm For $Y\subseteq X$, let $Y_{S}$ be the set of all subsets
\begin{equation*}
\{a_{i},a_{i+1},\ldots ,a_{i+r}\}\text{ }(i,r\in \{1,\ldots ,n\})
\end{equation*}%
of $Y$ such that $a_{i-1}\notin Y$ $($or $i=1)$ and $a_{i+r+1}\notin Y$ $($%
or $i+r=n)$.
\end{definition}

\begin{definition}
\rm Let $\alpha \in IF_n$ and let $k\in \mathbb{N}$. Then we put
$$\alpha (k) :=\{A\in (\dom \alpha )_{S}:\left\vert A\right\vert =k\},$$
$$\alpha^{o}(2k+1) :=\{\{a_{i},\ldots ,a_{i+2k}\}\in \alpha (2k+1):i\in 2\mathbb{N}-1\}.$$
\end{definition}
Note that $\alpha^{o}(2k+1)\subseteq \alpha (2k+1)$. \newline

For a set $M$ of natural numbers, let $\max M$ (let $\min M$) be the
greatest (the least) natural number in $M$ with respect to the natural order
in $\mathbb{N}$.

\begin{proposition}\label{prop9}
\rm Let $\alpha, \beta \in IF_n$. Then the following statements are equivalent:
\newline
(i) $\alpha \mathcal{J}\beta $. \newline
(ii) $\left\vert \alpha (k)\right\vert =\left\vert \beta (k)\right\vert $ and $%
\left\vert \alpha ^{o}(2k+1)\right\vert =\left\vert \beta^{o}(2k+1)\right\vert $ for all $k\in \mathbb{N}$.
\end{proposition}
\Pr
Suppose that $\alpha \mathcal{J}\beta $. Then there are $\gamma ,\delta,
\gamma _{1},\delta _{1},\in IF_n$ such that $\beta =\gamma \alpha \delta$
and $\alpha =\gamma _{1}\beta \delta _{1}$.
We have $\rank \alpha = \rank \beta$ since $IF_n \leq I_n$. Then from $\alpha = \gamma_1\beta\delta_1$ and $\beta = \gamma\alpha\delta$,
we obtain $|(\dom \alpha )_{S}| = |(\dom \beta)_{S}|$, and in particular, $\left\vert \alpha (k)\right\vert =\left\vert \beta (k)\right\vert $
for all $k \in \mathbb{N}$. Moreover, if $k \in \mathbb{N}$ and $B \in \beta(k)$ then we observe $B\gamma \in (\dom \alpha )_{S}$ and thus
$B\gamma \in \alpha (k)$.

Let $k\in 2\mathbb{N}+1$ and $B:=\{a_{i},\ldots ,a_{i+k-1}\}\in \beta
^{o}(k) $ (for some $i\in 2\mathbb{N}-1$). We have $B\gamma \in \alpha (k)$
and we will show that $B\gamma \in \alpha^o(k)$. Since $i$ is odd, we have $%
a_{i}<_{f}a_{i+1}>_{f}\cdots <_{f}a_{i+k-2}>_{f}a_{i+k-1}$. This implies $%
a_{i}\gamma <_{f}a_{i+1}\gamma >_{f}\cdots <_{f}a_{i+k-2}\gamma
>_{f}a_{i+k-1}\gamma $ and there is $l\in \{1,\ldots ,n\}$ with $\ $either $%
a_{i}\gamma =a_{l}$ and $a_{i+k-1}\gamma =a_{l+k-1}$ or $a_{i}\gamma
=a_{l+k-1}$ and $a_{i+k-1}\gamma =a_{l}$. This gives $a_{l}<_{f}a_{l+1}$ and
$l\in 2\mathbb{N}-1$, and consequently, $B\gamma \in \alpha ^{o}(k)$. This
shows $\left\vert \beta ^{o}(k)\right\vert \leq \left\vert \alpha
^{o}(k)\right\vert $. Dually, we can verify the converse inequation. Thus, $%
\left\vert \alpha ^{o}(k)\right\vert =\left\vert \beta ^{o}(k)\right\vert $.%
\newline

Conversely, let $\left\vert \alpha (k)\right\vert =\left\vert \beta
(k)\right\vert $ and $\left\vert \alpha ^{o}(2k+1)\right\vert =\left\vert
\beta ^{o}(2k+1)\right\vert $ for all $k\in \mathbb{N}$. Then for all $k\in
\mathbb{N}$, there is a bijection $f_{k}:\beta (k)\rightarrow \alpha (k)$
such that $f_{2k+1}(B)\in \alpha ^{o}(2k+1)$ for all $B\in \beta ^{o}(2k+1)$%
. We define now a mapping $\gamma :\dom \beta \rightarrow \dom \alpha $. For $%
k\in \mathbb{N}$, $B=\{a_{i},\ldots ,a_{i+k-1}\}\in \beta (k)$, and
$f_k(B)=\{a_{l},\ldots ,a_{l+k-1}\}$ (with $i,l \in \{1,\ldots ,n\}$) let
\begin{equation*}
a_{i+r}\gamma :=\left\{
\begin{array}{ll}
a_{l+r} & \text{if }k=1\text{ or }i\text{ and }l\text{ have the same parity}
\\
a_{l+k-(r+1)} & \text{otherwise}%
\end{array}%
\right.
\end{equation*}%
for $0\leq r\leq k-1$. The mapping $\gamma $ is well defined since $%
\dom \beta =\overset{p}{\underset{j=1}{\bigcup }}\beta (j)$, where $p:=\max
\{k\in \mathbb{N}:\beta (k)\neq \emptyset \}$.

We have to show that $\gamma \in IF_n$. For this let again $%
B=\{a_{i},\ldots ,a_{i+k-1}\}\in \beta (k)$ and $\{a_{l},\ldots
,a_{l+k-1}\}=f_{k}(B)$ for some $i,l,k\in \{1,\ldots ,n\}$.

We consider here the case $i\in 2\mathbb{N}-1$, the case $i\in 2\mathbb{N}$ can be handled in the same matter.
Suppose that $k\in 2\mathbb{N}+1$. Then $B\in \beta
^{o}(k)$ and $f_{k}(B)\in \alpha ^{o}(k)$, i.e. $l\in 2\mathbb{N}-1$. \
Since $i$ and $l$ are odd, we have $a_{i}<_{f}a_{i+1}>_{f}\cdots
<_{f}a_{i+k-2}>_{f}a_{i+k-1}$ and $a_{l}<_{f}a_{l+1}>_{f}\cdots
<_{f}a_{l+k-2}>_{f}a_{l+k-1}$, i.e. $a_{i}\gamma <_{f}a_{i+1}\gamma
>_{f}\cdots <_{f}a_{i+k-2}\gamma >_{f}a_{i+k-1}\gamma $.

Now, suppose that $k\in 2\mathbb{N}$. Since $i$
is odd, we have $a_{i}<_{f}a_{i+1}>_{f}\cdots >_{f}a_{i+k-2}<_{f}a_{i+k-1}$.
If $l$ is odd, then $a_{l}<_{f}a_{l+1}>_{f}\cdots
>_{f}a_{l+k-2}<_{f}a_{l+k-1}$, i.e. $a_{i}\gamma <_{f}a_{i+1}\gamma
>_{f}\cdots >_{f}a_{i+k-2}\gamma <_{f}a_{i+k-1}\gamma $. If $l$ is even then
$a_{l}>_{f}a_{l+1}<_{f}\cdots <_{f}a_{l+k-2}>_{f}a_{l+k-1}$, i.e. $%
a_{i+k-1}\gamma >_{f}a_{i+k-2}\gamma <_{f}\cdots <_{f}a_{i+1}\gamma
>_{f}a_{i}\gamma $.

This shows that $\gamma \in PFI_n$. Let $r\in \{1,\ldots ,n-1\}$ with $%
a_{r},a_{r+1}\in A$ for some $A\in (\im \gamma )_{S}$. We observe that $%
\{B\gamma :B\in (\dom \beta )_{S}\}=(\dom \alpha )_{S}$. Thus, there is $B\in
(\dom \beta )_{S}$ such that $B\gamma =A$ and there is $s\in \{1,\ldots ,n\}$
with $a_{s}=a_{r}\gamma ^{-1}$. Then $a_{r+1\text{ }}\gamma ^{-1}\in
\{a_{s+1},a_{s-1}\}$. If $r$ is odd then $a_{r}<_{f}a_{r+1}$. Assume that $%
a_{r}\gamma ^{-1}>_{f}a_{r+1}\gamma ^{-1}$. Then $s$ is even, i.e. $%
a_{s}\gamma >_{f}a_{s+1}\gamma $ (if $a_{r+1\text{ }}\gamma ^{-1}=$ $a_{s+1}$%
) and $a_{s-1}\gamma <_{f}a_{s}\gamma $ (if $a_{r+1\text{ }}\gamma ^{-1}=$ $%
a_{s-1}$). This gives $a_{r}>_{f}a_{r+1}$, a contradiction. If $r$ is even
then $a_{r}>_{f}a_{r+1}$ and we obtain $a_{r}\gamma ^{-1}>_{f}a_{r+1}\gamma
^{-1}$ by the same arguments. This provides $\gamma ^{-1}\in PFI_n$, i.e. $%
\gamma \in IF_n$.

Finally, we define $\delta :\im \alpha \rightarrow \im \beta $ by
\begin{equation*}
\delta :=\alpha ^{-1}\gamma ^{-1}\beta.
\end{equation*}%
Since $\alpha ,\beta ,\gamma \in IF_n$, we have $%
\delta =\alpha ^{-1}\gamma ^{-1}\beta \in IF_n$.

There holds $\beta =\gamma \alpha \delta $. In fact, for $a\in \dom \beta $,
we obtain $a\gamma \alpha \delta =a\gamma \alpha \alpha ^{-1}\gamma
^{-1}\beta =a\beta $ since $\dom \alpha =\im \gamma $ and $\dom \gamma =\dom \beta $.
\fbx

\section{Generating sets}

For convenience, we arrange such that $X_n$ is the set of the first positive integers $n$ for some $n \in \mathbb{N}$,
i.e. $X_n = \{1, \ldots, n\}$ with
$$1 <_f 2 >_f 3 <_f \cdots n.$$
Clearly, the minimal elements of the fence $\textbf{F} = (X_n, <_f)$ are odd and maximal elements are even.
For $a, b \in X_n$, we will write $a \equiv b~ (\mod 2)$ or shorter $a \equiv_2 b$ if $a$ and $b$ have the same parity.
Further, we denote by $\varepsilon_i$ the identity mapping on $X_n\setminus\{i\}$ for $i = 1, \ldots, n$,
i.e. $\varepsilon_i := id|_{X_n\setminus \{i\}}$.
\begin{notation}\label{not4}
\rm Let us put
$$J := \{\alpha \in IF_n : \rank \alpha \geq n-2\}.$$
\end{notation}

The aim of this section is to show that $J$ is a generating set for the semigroup $IF_n$.
Note, $\varepsilon_i^{-1} = \varepsilon_i \in J$ for $1 \leq i \leq n$.
\begin{lemma}\label{le6}
\rm Let $m,p\in \mathbb{N}$ with $m+p\leq n$ and $m \equiv_2 m+p$. Then
\begin{equation*}
\alpha =\left(
\begin{array}{ccccccccc}
1 & \cdots & m-2 & m & \cdots & m+p & m+p+2 & \cdots & n \\
1 & \cdots & m-2 & m+p & \cdots & m & m+p+2 & \cdots & n%
\end{array}%
\right) \in J \text{ }
\end{equation*}%
and $\alpha ^{-1}\in J$.
\end{lemma}
\Pr
By simple calculations, one can see that $\alpha \in IF_n$. Since $\rank \alpha = n-2$ and
$\alpha ^{-1}=\alpha$, we obtain $\alpha, \alpha ^{-1} \in J$.
\fbx

\begin{lemma}\label{le2}
\rm Let $m,p\in \mathbb{N}$ such that $m+p+2\leq n$. Then
\begin{equation*}
\alpha =\left(
\begin{array}{ccccccccc}
1 & \cdots & m-2 & m & \cdots & m+p & m+p+4 & \cdots & n \\
1 & \cdots & m-2 & m+2 & \cdots & m+p+2 & m+p+4 & \cdots & n%
\end{array}%
\right) \in \left\langle J \right\rangle \text{ }
\end{equation*}%
and $\alpha ^{-1}\in \left\langle J \right\rangle $.
\end{lemma}
\Pr
We have to consider two cases.

1) Suppose that $p$ is even. Then $m \equiv_2 m+p$ and we consider the following
transformations with $\rank \geq n-2$:
$$\beta _{1} =\left(
\begin{array}{ccccccccc}
1 & \cdots & m-2 & m & \cdots & m+p+2 & m+p+4 & \cdots & n \\
1 & \cdots & m-2 & m+p+2 & \cdots & m & m+p+4 & \cdots & n%
\end{array}%
\right)$$
and
$$\beta _{2} =\left(
\begin{array}{ccccccccc}
1 & \cdots & m & m+2 & \cdots & m+p+2 & m+p+4 & \cdots & n \\
1 & \cdots & m & m+p+2 & \cdots & m+2 & m+p+4 & \cdots & n%
\end{array}%
\right).$$
Clearly, $\beta_1, \beta_2 \in J$ by Lemma \ref{le6} and it is easy to verify that
\begin{equation*}
\alpha =\beta _{1}\beta _{2}\varepsilon_m \text{ and } \alpha
^{-1}=\varepsilon_m\beta _{2}\beta _{1}
\end{equation*}%
where $\varepsilon_m \in J$. Thus, we obtain $\alpha ,\alpha ^{-1}\in
\left\langle J\right\rangle $.\newline

2) Now suppose that $p$ is odd. Then $m \not\equiv_2 m+p$ and we consider the following
transformations with $\rank \geq n-2$:
$$\beta _{3} =\left(
\begin{array}{ccccccccc}
1 & \cdots & m-2 & m & \cdots & m+p+1 & m+p+3 & \cdots & n \\
1 & \cdots & m-2 & m+p+1 & \cdots & m & m+p+3 & \cdots & n%
\end{array}%
\right)$$
and
$$\beta _{4} =\left(
\begin{array}{ccccccccc}
1 & \cdots & m-1 & m+1 & \cdots & m+p+2 & m+p+4 & \cdots & n \\
1 & \cdots & m-1 & m+p+2 & \cdots & m+1 & m+p+4 & \cdots & n%
\end{array}%
\right).$$
Clearly, $\beta_3, \beta_4 \in J$ by Lemma \ref{le6} and it is easy to verify that
\begin{equation*}
\alpha =\beta _{3}\beta _{4} \text{ and } \alpha ^{-1}=\beta _{4}\beta_{3}.
\end{equation*}%
Thus, $\alpha ,\alpha ^{-1}\in \left\langle J\right\rangle $.
\fbx
\begin{corollary}\label{cor2}
\rm Let $m,p,k\in \mathbb{N}$ such that $m+p+2k\leq n$. Then
\begin{equation*}
\alpha =\left(
\begin{array}{ccccccccc}
1 & \cdots & m-2 & m & \cdots & m+p & m+p+2k+2 & \cdots & n \\
1 & \cdots & m-2 & m+2k & \cdots & m+p+2k & m+p+2k+2 & \cdots & n%
\end{array}%
\right) \in \left\langle J\right\rangle
\end{equation*}%
and $\alpha ^{-1}\in \left\langle J\right\rangle $.
\end{corollary}
\Pr
For $0\leq i<k$ we define the transformations
\begin{equation*}
\beta_{i}=\left(
\begin{array}{ccccccccc}
1 & \cdots & m+2i-2 & m+2i & \cdots & m+p+2i & m+p+2i+4 & \cdots & n \\
1 & \cdots & m+2i-2 & m+2i+2 & \cdots & m+p+2i+2 & m+p+2i+4 & \cdots & n%
\end{array}%
\right) \text{.}
\end{equation*}%
Note that $\beta _{i},\beta _{i}^{-1}\in \left\langle J\right\rangle $ ($%
0\leq i<k$) by Lemma \ref{le2}. It is easy to verify that $\alpha =\beta
_{0}\cdots \beta _{k-1}\in \left\langle J\right\rangle $ and $\alpha
^{-1}=\beta _{k-1}^{-1}\cdots \beta _{0}^{-1}\in \left\langle J\right\rangle$.
\fbx
\begin{lemma}\label{le3}
\rm Let $m,p\in \mathbb{N}$ such that $p$ is odd and $m+p+1\leq n$. Then
\begin{equation*}
\alpha =\left(
\begin{array}{ccccccccc}
1 & \cdots & m-2 & m & \cdots & m+p & m+p+3 & \cdots & n \\
1 & \cdots & m-2 & m+p+1 & \cdots & m+1 & m+p+3 & \cdots & n%
\end{array}%
\right) \in \left\langle J\right\rangle
\end{equation*}%
and $\alpha ^{-1}\in \left\langle J\right\rangle $.
\end{lemma}
\Pr
We define a transformation
\begin{equation*}
\beta_1 =\left(
\begin{array}{ccccccccc}
1 & \cdots & m-2 & m & \cdots & m+p+1 & m+p+3 & \cdots & n \\
1 & \cdots & m-2 & m+p+1 & \cdots & m & m+p+3 & \cdots & n%
\end{array}%
\right) \text{.}
\end{equation*}%
Clearly, $\beta_1 \in J$ by Lemma \ref{le6}. Then we can verify that $\alpha =\beta_1
\varepsilon_m \in \left\langle J\right\rangle $ and $\alpha
^{-1}=\varepsilon_m\beta_1\in \left\langle J\right\rangle $.
\fbx
\begin{corollary}\label{cor3}
\rm Let $m,p,k\in \mathbb{N}$ such that $p$ is odd and $m+p+2k-1\leq n$. Then
\begin{equation*}
\alpha =\left(
\begin{array}{ccccccccc}
1 & \cdots & m-2 & m & \cdots & m+p & m+p+2k+1 & \cdots & n \\
1 & \cdots & m-2 & m+p+2k-1 & \cdots & m+2k-1 & m+p+2k+1 & \cdots & n%
\end{array}%
\right) \in \left\langle J\right\rangle
\end{equation*}%
and $\alpha ^{-1}\in \left\langle J\right\rangle $.
\end{corollary}
\Pr
Let
\begin{equation*}
\beta _{1}=\left(
\begin{array}{ccccccccc}
1 & \cdots & m-2 & m & \cdots & m+p & m+p+2k & \cdots & n \\
1 & \cdots & m-2 & m+2k-2 & \cdots & m+p+2k-2 & m+p+2k & \cdots & n%
\end{array}%
\right)
\end{equation*}%
and
\begin{equation*}
\beta_{2}=\left(
\begin{array}{ccccccccc}
1 & \cdots & m+2k-4 & m+2k-2 & \cdots & m+p+2k-2 & m+p+2k+1 & \cdots & n \\
1 & \cdots & m+2k-4 & m+p+2k-1 & \cdots & m+2k-1 & m+p+2k+1 & \cdots & n%
\end{array}%
\right).
\end{equation*}%
Note that $\beta _{1}\in \left\langle J\right\rangle $ (by Corollary \ref{cor2}) and
$\beta _{2}\in \left\langle J\right\rangle $ (by Lemma \ref{le3}). It is easy to
verify that $\alpha =\beta _{1}\beta _{2}$ and $\alpha ^{-1}=\beta
_{2}^{-1}\beta_{1}^{-1}$, and thus $\alpha ,\alpha ^{-1}\in \left\langle
J\right\rangle$.
\fbx
\begin{lemma}\label{le4}
\rm Let $m,p,k\in \mathbb{N}$ such that $p$ is even and $m+p+2k\leq n$. Then
\begin{equation*}
\alpha =\left(
\begin{array}{ccccccccc}
1 & \cdots & m-2 & m & \cdots & m+p & m+p+2k+2 & \cdots & n \\
1 & \cdots & m-2 & m+p+2k & \cdots & m+2k & m+p+2k+2 & \cdots & n%
\end{array}%
\right) \in \left\langle J\right\rangle
\end{equation*}%
and $\alpha ^{-1}\in \left\langle J\right\rangle $.
\end{lemma}
\Pr
Let
\begin{equation*}
\beta_1=\left(
\begin{array}{ccccccccc}
1 & \cdots & m-2 & m & \cdots & m+p & m+p+2k+2 & \cdots & n \\
1 & \cdots & m-2 & m+2k & \cdots & m+p+2k & m+p+2k+2 & \cdots & n%
\end{array}%
\right)
\end{equation*}%
and
\begin{equation*}
\beta_2=\left(
\begin{array}{ccccccccc}
1 & \cdots & m+2k-2 & m+2k & \cdots & m+p+2k & m+p+2k+2 & \cdots & n \\
1 & \cdots & m+2k-2 & m+p+2k & \cdots & m+2k & m+p+2k+2 & \cdots & n%
\end{array}%
\right).
\end{equation*}
Note that $\beta _{1}\in \left\langle J\right\rangle $ (by Corollary \ref{cor2}) and
$\beta _{2}\in J$ (by Lemma \ref{le6}). It is easy to
verify that $\alpha =\beta _{1}\beta _{2}$ and $\alpha ^{-1}=\beta
_{2}\beta _{1}^{-1}$, and thus $\alpha ,\alpha ^{-1}\in \left\langle
J\right\rangle $. \fbx
\begin{lemma}\label{le5}
\rm Let $Y\subseteq X_n$. Then $id|_{X_n\setminus Y}\in \left\langle J\right\rangle$.
\end{lemma}
\Pr
If $Y=\emptyset$, i.e. $X_n\setminus Y = X_n$, then $id|_{X_n} = id \in J$.
Let $\emptyset \neq Y:=\{i_{1},\ldots,i_{k}\}\subseteq X_n$ with $k \in \{1,\ldots ,n\}$. Then it is easy
to verify that $id|_{X_n\setminus Y}=\varepsilon_{i_1}\cdots \varepsilon_{i_k} \in \left\langle J\right\rangle$.
\fbx

\begin{proposition}\label{prop1}
\rm Let $\alpha \in IF_n$.
Then there are transformations $\eta_1, \ldots, \eta_k$, $\eta_{k+1}, \ldots, \eta_l \in J$ ($k<l \in \mathbb{N}$) such that
$\eta_1^{-1}, \ldots, \eta_k^{-1}, \eta_{k+1}^{-1}, \ldots, \eta_l^{-1} \in J$, $\dom\alpha \subseteq \im (\eta_1 \cdots \eta_k)$,
$\im\alpha \subseteq \dom (\eta_{k+1} \cdots \eta_l)$, and $x(\eta_1\ldots\eta_k\alpha\eta_{k+1}\ldots\eta_l) \equiv_2 x$ for all
$x \in \dom (\eta_1\ldots\eta_k\alpha\eta_{k+1}\ldots\eta_l)$.
\end{proposition}
\Pr If $a \equiv_2 a\alpha$ for all $a \in \dom \alpha$ then $id|_{\dom\alpha}\alpha~ id|_{\im\alpha} = \alpha$.
This shows the assertion, since $id|_{\dom\alpha}, id|_{\im\alpha} \in \langle J \rangle$ by Lemma \ref{le5}.

Let $a \in \dom \alpha$ such that $a \not\equiv_2 a\alpha$. Then it is clear that $a-1, a+1 \notin \dom \alpha$.

If $a$ is even then we put
$$\eta = \left(
             \begin{array}{ccccccc}
               1 & 3 & \cdots & a & a+2 & \cdots & n \\
               a & 1 & \cdots & a-2 & a+2 & \cdots & n \\
             \end{array}
           \right) \in J.$$
We observe that $\eta^{-1} \in \langle J \rangle$. Moreover, it is easy to see that $\im \alpha = \im(\eta\alpha)$, $x\alpha^{-1} \equiv_2 x(\eta\alpha)^{-1}$
for all $x \in \im\alpha\setminus\{a\alpha\}$ and $a\alpha\alpha^{-1} = a \not\equiv_2 1 = a\eta^{-1} = a\alpha(\eta\alpha)^{-1}$.
This shows that
$$|\{x \in \im\alpha : x \not\equiv_2 x(\eta\alpha)^{-1}\}| = |\{x \in \im\alpha : x \not\equiv_2 x\alpha^{-1}\}| -1.$$

If $a$ is odd then $a\alpha$ is even and we put
$$\eta = \left(
             \begin{array}{cccccccc}
               1 & \cdots & a\alpha-2 & a\alpha & a\alpha+2 & \cdots & n \\
               3 & \cdots & a\alpha & 1 & a\alpha+2 & \cdots & n \\
             \end{array}
           \right) \in J,$$
with $\eta^{-1} \in \langle J \rangle$. By dual arguments we obtain
$$|\{x \in \dom\alpha : x \not\equiv_2 x(\alpha\eta)\}| = |\{x \in \dom\alpha : x \not\equiv_2 x\alpha\}| -1.$$

Continuing in this way, starting with the even cases, we obtain transformations $\eta_1, \ldots, \eta_k$, $\eta_{k+1}, \ldots, \eta_l \in J$ ($k<l \in \mathbb{N}$) such that
$\eta_1^{-1}, \ldots, \eta_k^{-1}$, $\eta_{k+1}^{-1}, \ldots, \eta_l^{-1} \in J$ and $x(\eta_1\ldots\eta_k\alpha\eta_{k+1}\ldots\eta_l) \equiv_2 x$ for all
$x \in \dom (\eta_1 \ldots \eta_k \alpha \eta_{k+1} \ldots \eta_l)$.\\
\fbx
\begin{notation}\label{not1}
\rm Let $\alpha \in PFI_n$ and let $A,B\in (\dom \alpha )_{S}$ (or $A,B\in (\im\alpha)_{S}$). Then we write $A < B$ if all elements in $A$ are less than
any element in $B$ with respect to the natural order of$\ \mathbb{N}$. Further,
we write
\begin{equation*}
A\prec B
\end{equation*}%
if $A < B$ and for each $C\in (\dom\alpha)_{S}$ (for each $C\in
(\im\alpha)_{S}$, respectively) the following implication holds: $A \leq C \leq B\Rightarrow A=C$ or $B=C$.
\end{notation}

Any transformation $\alpha \in IF_n$ with $a \equiv_2 a\alpha$ for all $a \in \dom\alpha$
can be written in the following form:
\begin{notation}\label{not2}
\rm Let
$$\alpha = \left(
                \begin{array}{ccccccc}
                  A_1 & \prec \cdots \prec & A_{i-1} & \prec & A_i & \prec \cdots \prec & A_p \\
                  A_1 & \prec \cdots \prec & A_{i-1} & < & B_i & \cdots & B_p \\
                \end{array}
              \right) \in IF_n,$$
with $i \leq p \in \{1, \ldots, n\}$, and $a \equiv_2 a\alpha$ for all $a \in \dom \alpha$ such that $i=1$ or \\
(i) $a\alpha = a$ for all $a \in A_1 \cup \cdots \cup A_{i-1}$ and\\
(ii) $A_{i-1} < B_l$ for all $l \in \{i, \ldots, p\}$.\\
\\
Further, let
$$r_{j}:=\min A_{j},~~ s_{j}:=\max A_{j},~~ t_{j}:=\min B_{j},~~ u_{j}:=\max B_{j},$$
for $1\leq j\leq p$.
\end{notation}
\begin{proposition}\label{prop2}
\rm Let $\alpha$ be as in Notation \ref{not2}. Then there exist $\omega_1, \omega_2 \in \langle J \rangle$
with $\omega_1^{-1}, \omega_2^{-1} \in \langle J \rangle$, $\dom \alpha \subseteq \im \omega_1$, $\im \alpha \subseteq \dom \omega_2$ such that $\omega_1\alpha\omega_2$ has the form
$$\omega_1\alpha\omega_2 = \left(
                \begin{array}{ccccccccc}
                  A_1 & \prec \cdots \prec & A_{i-1} & \prec & A_i' & \prec & A_{i+1}' & \prec \cdots \prec & A_p' \\
                  A_1 & \prec \cdots \prec & A_{i-1} & \prec & B_i' & < & B_{i+1}' & \cdots & B_p' \\
                \end{array}
              \right) \in IF_n,$$
with $a \equiv_2 a(\omega_1\alpha\omega_2)$ for all $a \in \dom (\omega_1\alpha\omega_2)$, and $B_{i}' < B_l'$ for all $l \in \{i+1, \ldots, p\}$ such that $i=1$ or $a(\omega_1\alpha\omega_2) = a$ for all $a \in A_1 \cup \cdots \cup A_{i-1}$.
\end{proposition}
\Pr
We will define the transformations $\omega_1$ and $\omega_2$ with $\dom \alpha \subseteq \im \omega_1$ and $\im \alpha \subseteq \dom \omega_2$
such that $\omega_1\alpha\omega_2$ is the required mapping of our assertion.
The concrete calculations we leave to the reader. \\
Let $k \in \{i, \ldots, p\}$ such that $A_{i-1} \prec B_k$ if $i > 1$, and $B_{k} < B_l$ for all $l \in \{1, \ldots, p\}\setminus\{k\}$ if $i=1$, respectively.
Note that if $k=i$ then $\omega_1 = \omega_2 = id \in J$. Thus, let $k>i$. Then we consider the following seven cases.
Note that the cases are not mutually exclusive (i.e. the transformation $\alpha$ can satisfy more than one case), but cover all the possibilities. \\
\\
1. If $r_i \equiv_2 s_k$ then we put $\omega_2 = id$ and
$$\omega_1 = \left(
             \begin{array}{ccccccccc}
               1 & \cdots & r_i-2 & r_i & \cdots & s_k & s_k+2 & \cdots & n \\
               1 & \cdots & r_i-2 & s_k & \cdots & r_i & s_k+2 & \cdots & n \\
             \end{array}
           \right).$$
2. If $r_i \not\equiv_2 s_k$ and $r_i-2 \notin \dom \alpha$ (or $r_i-1 = 1$) then we put $\omega_2 = id$ and
$$\omega_1 = \left(
             \begin{array}{ccccccccc}
               1 & \cdots & r_i-3 & r_i-1 & \cdots & s_k & s_k+2 & \cdots & n \\
               1 & \cdots & r_i-3 & s_k & \cdots & r_i-1 & s_k+2 & \cdots & n \\
             \end{array}
           \right).$$
3. If $r_i \not\equiv_2 s_k$ and $s_k+2 \notin \dom \alpha$ (or $s_k+1 = n$) then we put $\omega_2 = id$ and
$$\omega_1 = \left(
             \begin{array}{ccccccccc}
               1 & \cdots & r_i-2 & r_i & \cdots & s_k+1 & s_k+3 & \cdots & n \\
               1 & \cdots & r_i-2 & s_k+1 & \cdots & r_i & s_k+3 & \cdots & n \\
             \end{array}
           \right).$$
4. If $u_i \equiv_2 t_k$ then we put $\omega_1 = id$ and
$$\omega_2 = \left(
             \begin{array}{ccccccccc}
               1 & \cdots & t_k-2 & t_k & \cdots & u_i & u_i+2 & \cdots & n \\
               1 & \cdots & t_k-2 & u_i & \cdots & t_k & u_i+2 & \cdots & n \\
             \end{array}
           \right).$$
5. If $u_i \not\equiv_2 t_k$ and $t_k-2 \notin \im \alpha$ (or $t_k-1 = 1$) then we put
$\omega_1 = id$ and
$$\omega_2 = \left(
             \begin{array}{ccccccccc}
               1 & \cdots & t_k-3 & t_k-1 & \cdots & u_i & u_i+2 & \cdots & n \\
               1 & \cdots & t_k-3 & u_i & \cdots & t_k-1 & u_i+2 & \cdots & n \\
             \end{array}
           \right).$$
6. If $u_i \not\equiv_2 t_k$ and $u_i+2 \notin \im \alpha$ (or $u_i+1 = n$) then we put
$\omega_1 = id$ and
$$\omega_2 = \left(
             \begin{array}{ccccccccc}
               1 & \cdots & t_k-2 & t_k & \cdots & u_i+1 & u_i+3 & \cdots & n \\
               1 & \cdots & t_k-2 & u_i+1 & \cdots & t_k & u_i+3 & \cdots & n \\
             \end{array}
           \right).$$

Clearly, $\dom \alpha \subseteq \im \omega_1$, $\im \alpha \subseteq \dom \omega_2$, and $\omega_1, \omega_2 \in J$ (by Lemma \ref{le6})
for all cases 1. - 6.\\
\\
7. It remains the case $r_i \not\equiv_2 s_k$ and $u_i \not\equiv_2 t_k$ and $r_i-2, s_k+2 \in \dom \alpha$ and $t_k-2, u_i+2 \in \im \alpha$
where $1 = r_1 \in \dom \alpha$ and $1 = t_k \in \im \alpha$ in the case $i=1$.

7.1. Let $k = i+1$. First, we will show that $r_{i} = t_{i+1}$. In the case $i=1$, it is clear.
For the case $i>1$, we have that $t_{i+1} = u_{i-1}+2$ (since $A_{i-1} \prec B_{i+1}$ and $t_{i+1}-2 \in \im \alpha$),
$r_{i} = s_{i-1}+2$ (since $A_{i-1} < A_{i}$ and $r_{i}-2 \in \dom \alpha$) and $u_{i-1}+2 = s_{i-1}+2$ (since $a\alpha = a$ for all $a \in A_1 \cup \cdots \cup A_{i-1}$).
Altogether, we obtain $r_{i} = t_{i+1}$. Since $r_i \not\equiv_2 s_{i+1}$, we have $r_{i} = t_{i+1}  = r_{i+1}\alpha \equiv_2 r_{i+1}$. Thus, we get $r_{i+1} \not\equiv_2 s_{i+1}$ and
we put $\omega_1 = \eta_1\eta_2$ and $\omega_2 = id$, where
$$\eta_1 = \left(
             \begin{array}{ccccccccc}
               1 & \cdots & r_i-2 & r_i & \cdots & s_{i+1}-1 & s_{i+1}+1 & \cdots & n \\
               1 & \cdots & r_i-2 & s_{i+1}-1 & \cdots & r_i & s_{i+1}+1 & \cdots & n \\
             \end{array}
           \right)$$
$$\eta_2 = \left(
             \begin{array}{ccccccccc}
               1 & \cdots & r_{i+1}-3 & r_{i+1}-1 & \cdots & s_{i+1}-1 & s_{i+1}+2 & \cdots & n \\
               1 & \cdots & r_{i+1}-3 & s_{i+1} & \cdots & r_{i+1} & s_{i+1}+2 & \cdots & n \\
             \end{array}
           \right).$$
Clearly, $\eta_1 \in J$ by Lemma \ref{le6} and $\eta_2 \in \langle J \rangle$ by Lemma \ref{le3}.
Note that $r_{i+1}-2 \notin \dom \alpha$, since otherwise
$s_{i}= r_{i+1}-2 \equiv_2 r_{i+1} \equiv_2 r_{i}$ implies $u_{i} \equiv_2 r_i = t_{i+1}$ which is a contradiction.
Thus, it is easy to verify that $\dom \alpha \subseteq \im \omega_1$.

7.2. Let $k > i+1$. We define a transformation $\tau$ as following:

a) If $r_{i+1} \equiv_2 s_k$ then we put
$$\tau = \left(
             \begin{array}{ccccccccc}
               1 & \cdots & r_{i+1}-2 & r_{i+1} & \cdots & s_k & s_k+2 & \cdots & n \\
               1 & \cdots & r_{i+1}-2 & s_k & \cdots & r_{i+1} & s_k+2 & \cdots & n \\
             \end{array}
           \right).$$

b) If $r_{i+1} \not\equiv_2 s_k$, i.e. $r_i \equiv_2 r_{i+1}$ then we put
$$\tau = \left(
             \begin{array}{ccccccccc}
               1 & \cdots & r_{i+1}-3 & r_{i+1}-1 & \cdots & s_k & s_k+2 & \cdots & n \\
               1 & \cdots & r_{i+1}-3 & s_k & \cdots & r_{i+1}-1 & s_k+2 & \cdots & n \\
             \end{array}
           \right)$$
By Lemma \ref{le6}, we have $\tau, \tau^{-1} \in J$. We have to verify that $r_{i+1}-2 \notin \dom \alpha$.
Assume the opposite that $r_{i+1}-2 \in \dom \alpha$. Then $s_i = r_{i+1}-2$ and thus $s_i = r_{i+1}-2 \equiv_2 r_{i+1}\equiv_2 r_i$.
Therefore, we have $r_i \equiv_2 s_i \equiv_2 t_i \equiv_2 u_i$.
Moreover, we have $r_i = t_k = 1$ in the case $i=1$. If $i>1$ then $u_{i-1} = s_{i-1} = r_i-2 \equiv_2 r_i$ and $u_{i-1} \equiv_2 u_{i-1}+2 = t_{k}$ implies $r_i \equiv_2 t_{k}$. Thus, we obtain $u_i \equiv_2 t_{k}$, a contradiction. Hence, $\dom \alpha \subseteq \im \tau$.\\

Now, we consider the transformation
$$\tau\alpha = \left(
                \begin{array}{cccccccccc}
                  A_1 & \cdots & A_{i-1} & A_i & A_{i+1}^* & \cdots & A_k^* & A_{k+1} & \cdots & A_p \\
                  A_1 & \cdots & A_{i-1} & B_i & B_{i+1}^* & \cdots & B_k^* & B_{k+1} & \cdots & B_p \\
                \end{array}
              \right) \in IF_n,$$
with $A_{i-1} \prec B_{i+1}^*$. For this transformation, we have the case 7.1. with corresponding transformations $\eta_1, \eta_2 \in \langle J\rangle$.
Then we put $\omega_1 = \eta_1\eta_2\tau$ and $\omega_2 = id$ with $\omega_1^{-1}, \omega_2^{-1} \in \langle J\rangle$, $\dom \alpha \subseteq \im \omega_1$
and $\im \alpha \subseteq \dom \omega_2$.\\
\fbx
\begin{proposition}\label{prop3}
\rm Let $\alpha$ be as in Notation \ref{not2} with $A_{i-1} \prec B_i$. Then there exist $\omega_1, \omega_2 \in \langle J \rangle$ such that
$\omega_1^{-1}, \omega_2^{-1} \in \langle J\rangle$, $\dom \alpha \subseteq \im \omega_1$, $\im \alpha \subseteq \dom \omega_2$, and
$$\omega_1\alpha\omega_2 = \left(
                \begin{array}{ccccccccc}
                  A_1 & \prec \cdots \prec & A_{i-1} & \prec & A_i' & \prec & A_{i+1}' & \prec \cdots \prec & A_p' \\
                  A_1 & \prec \cdots \prec & A_{i-1} & \prec & A_i' & < & B_{i+1}' & \cdots & B_p' \\
                \end{array}
              \right) \in IF_n,$$
with $a(\omega_1\alpha\omega_2) = a$ for all $a \in (A_1 \cup \cdots \cup A_{i-1} \cup A_i')$.
\end{proposition}
\Pr If $a\alpha = a$ for all $a \in A_i$ then $\omega_1 = \omega_2 =id$. Let $a\alpha\neq a$ for some $a \in A_i$. Then we put
$$\eta_1 = \left(
             \begin{array}{ccccccccc}
               1 & \cdots & t_{i}-2 & r_{i}\alpha & \cdots & s_i\alpha & s_i+2 & \cdots & n \\
               1 & \cdots & t_{i}-2 & r_i & \cdots & s_{i} & s_i+2 & \cdots & n \\
             \end{array}
           \right), ~~\mbox{ if } r_i \geq t_i \mbox{ and}$$
$$\eta_2 = \left(
             \begin{array}{ccccccccc}
               1 & \cdots & r_{i}-2 & r_{i}\alpha & \cdots & s_i\alpha & u_i+2 & \cdots & n \\
               1 & \cdots & r_{i}-2 & r_i & \cdots & s_{i} & u_i+2 & \cdots & n \\
             \end{array}
           \right), ~~\mbox{ if } r_i \leq t_i.$$
Clearly, $\eta_1,\eta_2,\eta_1^{-1},\eta_2^{-1} \in \langle J \rangle$ by Corollary \ref{cor2} (if $r_i\alpha=t_i$) or Corollary \ref{cor3} (if $r_i\alpha=u_i$ and $r_i \not\equiv_2 s_i$) or Lemma \ref{le4} (if $r_i\alpha=u_i$ and $r_i \equiv_2 s_i$).
If $r_i \geq t_i$ then $\dom \alpha \subseteq \im \eta_1$ and we put $\omega_1 = \eta_1$ and $\omega_2 = id$. If $r_i \leq t_i$ then $\im \alpha \subseteq \dom \eta_2$ and we put $\omega_1 = id$ and $\omega_2 = \eta_2$.
\fbx\\

From Proposition \ref{prop1}, Proposition \ref{prop2}, and Proposition \ref{prop3} (frequently used) we obtain
\begin{corollary}\label{cor1}
\rm Let $\alpha \in IF_n$. Then there exist $\omega_1, \omega_2 \in \langle J \rangle$ such that $\omega_1^{-1}, \omega_2^{-1} \in \langle J \rangle$,
$\dom \alpha \subseteq \im \omega_1$, $\im \alpha \subseteq \dom \omega_2$,
and $a(\omega_1\alpha\omega_2) = a$ for all $a \in \dom(\omega_1\alpha\omega_2)$.
\end{corollary}
\begin{theorem}\label{theo1}
$IF_n = \langle J \rangle$.
\end{theorem}
\Pr Let $\alpha \in IF_n$. Then by Corollary \ref{cor1}, there exist $\omega_1, \omega_2 \in \langle J \rangle$ such that $\omega_1^{-1}, \omega_2^{-1} \in \langle J \rangle$, $\dom \alpha \subseteq \im \omega_1$, $\im \alpha \subseteq \dom \omega_2$,
and $a(\omega_1\alpha\omega_2) = a$ for all $a \in \dom(\omega_1\alpha\omega_2)$. Therefore, we have
$$\omega_1\alpha\omega_2 = \varepsilon_{i_1} \cdots \varepsilon_{i_k} \in \langle J \rangle$$
(by Lemma \ref{le5}), where $\{i_1, \ldots, i_k\} = X_n\setminus\dom(\omega_1\alpha\omega_2)$, $k \in \{1, \ldots, n\}$.

Finally, we obtain $\alpha \in \langle J \rangle$, since $\alpha = \omega_1^{-1}\omega_1\alpha\omega_2\omega_2^{-1}$.\\
\fbx

\section{Rank of the semigroup $IF_n$ for even $n$}

Let $n\in 2\mathbb{N}+1$. Using the GAP software, we have observed that $IF_n$ is not generated by the set
$\{\alpha \in IF_n : \rank \alpha \geq n-1\}$. Moreover, there is no least generating set for $IF_n$. But, in the case
$n$ is even the situation is different. There is a least generating set and all its elements have rank $\geq n-1$.

Throughout this section, let $n \in 2\mathbb{N}$ and let $X_n$ be again the up-fence $1 <_{f} 2 >_{f}\cdots <_{f} n$.
We describe the least generating set and calculate the rank of the semigroup $IF_n$.
\begin{notation}
\rm We put\\
\\

$\sigma _{1} :=\left(
\begin{array}{llll}
1 & 3 & \cdots  & n \\
n & 1 & \cdots  & n-2%
\end{array}%
\right)$; \\
\\

$\sigma_{2} := \sigma_{1}^{-1} = \left(
\begin{array}{llll}
1 & \cdots  & n-2 & n \\
3 & \cdots  & n & 1%
\end{array}%
\right)$; \\
\\

$\gamma _{i}:=\left(
\begin{array}{llllll}
1 & \cdots  & i-1 & i+1 & \cdots  & n \\
i-1 & \cdots  & 1 & i+1 & \cdots  & n%
\end{array}%
\right)$  for  $i\in 2\mathbb{N}, 4\leq i\leq n$;\\
\\

$\delta _{i}:=\left(
\begin{array}{llllll}
1 & \cdots  & i-1 & i+1 & \cdots  & n \\
1 & \cdots  & i-1 & n & \cdots  & i+1%
\end{array}%
\right)$  for  $i\in 2\mathbb{N}-1, 1\leq i\leq n-3$;\\
\\

$G := \{id\} \cup \{\sigma _{1}, \sigma_2\} \cup \{\gamma _{i} : i\in 2\mathbb{N}, 4\leq i\leq n\} \cup \{\delta _{i} : i\in 2\mathbb{N}-1, 1\leq i\leq n-3\}$.\\
\end{notation}

Note that $\sigma_2^{-1} = \sigma_1$,~~ $\gamma_i^{-1} = \gamma_i$, ~and~ $\delta_i^{-1} = \delta_i$.

\begin{theorem}\label{theo2}
\rm $IF_n = \langle G \rangle$.
\end{theorem}
\Pr From Theorem \ref{theo1}, we have $IF_n = \langle J \rangle$. It remains to show that $J \subseteq \langle G \rangle$.
For this we have to show that all transformations $\varepsilon_{i}$ for $i \in \{1, \ldots, n\}$ as well as all transformations which are used in Proposition \ref{prop1}, Proposition \ref{prop2}, and Proposition \ref{prop3} belong to $\langle G \rangle$.

We observe that $\varepsilon_{i}=\gamma _{i}\gamma _{i}$ for $i\in 2\mathbb{N}$, $4\leq i\leq n$ and
$\varepsilon_{i}=\delta _{i}\delta _{i}$ for $i\in 2\mathbb{N}-1$, $1\leq i\leq n-3$ as well as
$\varepsilon_{2}=\sigma_1\sigma_2$ and $\varepsilon_{n-1}=\sigma_2\sigma_1$.

For the transformations in Proposition \ref{prop1} we have
$$\eta = \left(
             \begin{array}{ccccccc}
               1 & 3 & \cdots & i & i+2 & \cdots & n \\
               i & 1 & \cdots & i-2 & i+2 & \cdots & n \\
             \end{array}
           \right) = \delta_{i+1}\sigma_1\delta_{i-1} \in \langle G \rangle,$$
if $i=a$ is even and $i < n$. If $i=n$ then $\eta = \sigma_1$. Further, we have
$$\eta = \left(
             \begin{array}{cccccccc}
               1 & \cdots & i-2 & i & i+2 & \cdots & n \\
               3 & \cdots & i & 1 & i+2 & \cdots & n \\
             \end{array}
           \right) =\delta_{i-1}\sigma_2\delta_{i+1} \in \langle G \rangle,$$
if $i=a\alpha$ is even.

For the transformations in Proposition \ref{prop2} we put
$$\beta_{i,j} := \left(
             \begin{array}{ccccccccc}
               1 & \cdots & i-1 & i+1 & \cdots & j-1 & j+1 & \cdots & n \\
               1 & \cdots & i-1 & j-1 & \cdots & i+1 & j+1 & \cdots & n \\
             \end{array}
           \right) = \beta_{i,j}^{-1},$$
for $1 \leq i < j \leq n$ and $i \equiv_2 j$. Clearly, $\beta_{i,j} \in \langle G \rangle$ since
$$\beta_{i,j} = \left\{\begin{array}{ll}
                  \delta_i\delta_{n-j+i+1}\delta_i, & \mbox{if } i \mbox{ and } j \mbox{ are odd}; \\
                  \gamma_j\gamma_{j-i}\gamma_j, & \mbox{if } i \mbox{ and } j \mbox{ are even}.
                \end{array}\right.$$
It is easy to verify that $\omega_1$, $\omega_2$, $\eta_1$ and $\tau$ are all of the form $\beta_{i,j}$ for suitable $i$ and $j$.

Further, we have
$$\eta_2 = \left(
             \begin{array}{ccccccccc}
               1 & \cdots & i-1 & i+1 & \cdots & j-2 & j+1 & \cdots & n \\
               1 & \cdots & i-1 & j-1 & \cdots & i+2 & j+1 & \cdots & n \\
             \end{array}
           \right),$$
for suitable $i$ and $j$, and $i+1 \not\equiv_2 j-2$. Clearly, $\eta_2 \in \langle G \rangle$ since $\eta_2 = \beta_{i,j}\varepsilon_{i+1}$.\\

For the transformations in Proposition \ref{prop3} we have
$\eta_1, \eta_2 \in \langle J \rangle$ by Corollary \ref{cor2} (if $r_i\alpha=t_i$) or Corollary \ref{cor3} (if $r_i\alpha=u_i$ and $r_i \not\equiv_2 s_i$) or Lemma \ref{le4} (if $r_i\alpha=u_i$ and $r_i \equiv_2 s_i$).\\

For the transformation $\alpha$ in Corollary \ref{cor2}, we have
\[\alpha =\left(
\begin{array}{lllllllll}
1 & \ldots  & i-1 & i+1 & \ldots  & j-1-2k & j+1 & \ldots  & n \\
1 & \ldots  & i-1 & i+1+2k & \ldots  & j-1 & j+1 & \ldots  & n%
\end{array}%
\right)
\]%
with $m=i+1$ and $m+p=j-1-2k$.
Hence, we obtain $\alpha \in \left\langle G \right\rangle $ since
$$\alpha = \left\{\begin{array}{ll}
                  \beta _{i,j}\beta_{i+2k,j}\varepsilon_{i+1}\ldots \varepsilon_{i+2k-1}, & \mbox{if } i \equiv_2 j; \\
                  \beta _{i,j-1}\beta_{i+2k-1,j}\varepsilon_{i+1}\ldots \varepsilon_{i+2k-2}, & \mbox{if } i \not\equiv_2 j.
                \end{array}\right.$$

For the transformation $\alpha$ in Corollary \ref{cor3}, we have
\[\alpha =\left(
\begin{array}{lllllllll}
1 & \ldots  & i-1 & i+1 & \ldots  & j-2k & j+1 & \ldots  & n \\
1 & \ldots  & i-1 & j-1 & \ldots  & i+2k & j+1 & \ldots  & n%
\end{array}%
\right)
\]%
with $m=i+1$ and $m+p=j-2k$. Then we can verify that $\alpha \in
\left\langle G\right\rangle$ since $\alpha =\beta _{i,j}\varepsilon_{i+1}\ldots\varepsilon_{i+2k-1}$. \\

For the transformation $\alpha$ in Lemma \ref{le4}, we have
\[\alpha =\left(
\begin{array}{lllllllll}
1 & \ldots  & i-1 & i+1 & \ldots  & j-2k-1 & j+1 & \ldots  & n \\
1 & \ldots  & i-1 & j-1 & \ldots  & i+2k+1 & j+1 & \ldots  & n%
\end{array}%
\right)
\]%
with $m=i+1$ and $m+p=j-2k-1$. We have $\alpha \in \left\langle
G\right\rangle $ since $\alpha =\beta _{i,j}\varepsilon_{i+1}\ldots \varepsilon_{i+2k}$.
\fbx
\begin{proposition}\label{prop4}
\rm The set $G$ is the least generating set for $IF_n$.
\end{proposition}
\Pr
Theorem \ref{theo2} shows that $G$ is a generating set for $IF_n$. Let $\alpha ,\beta \in G$ with $\alpha \neq \beta $ and $\{\alpha ,\beta
\}\neq \{\sigma_1 ,\sigma_2 \}$. It is easy to verify that $\rank \alpha \beta =n-2$.
Moreover, we observe $\rank \sigma_1^{2}=\rank \sigma_2^{2}=n-2$. Let $\alpha
=\alpha_{1}\cdots \alpha_{m}$ with $\alpha_{1},\ldots,\alpha_{m}\in IF_n$, $2\leq m\in \mathbb{N}$, such that $\rank \alpha
=n-1$. Without loss of generality, we can assume that $\alpha_{i}\neq
id$ for $1\leq i\leq m$. Then $\alpha_{1},\ldots,\alpha_{m}\in
\{\beta \in IF_n:\rank \beta =n-1\}$. Since $G$ is a generating set for $\{\beta \in IF_n:\rank \beta =n-1\}$,
there is $\rho \in G$ such that $\alpha_{1},\ldots,\alpha_{m}\in \{\rho^{j}:j\in \mathbb{N}\}$
or $\alpha_{1},\ldots,\alpha_{m}$ $\in \{\sigma_1,\sigma_2\}$ with $\alpha_{i}\neq \alpha_{i+1}$ for $1\leq i\leq m-1$. This shows
that any $\alpha \in G$ can not be generated by a set without
this $\alpha$. Thus, each generating set of $IF_n$ have to contain $G$ and the assertion is shown.
\fbx\\

Since $|G| = n+1$ from Theorem \ref{theo2} and Proposition \ref{prop4}, we obtain
\begin{theorem}\label{theo}
\rm Let $n\in 2\mathbb{N}$. Then $\rank IF_n = n+1$.
\end{theorem}

\textbf{Acknowledgements.}
The authors wish to thank the Alexander von Humboldt Foundation for the financial support and Assoc. Prof. Kalcho Todorov for his essential suggestions concerning this research.

\end{document}